\documentclass[a4paper,12pt]{article}
\usepackage{amssymb}
\title{Lagrangian Pairs\\ and\\ Lagrangian Orthogonal Matroids}
\date{8 September 2002}
\author{Richard F. Booth\\
{\small         Department of Mathematics,}\\
{\small         UMIST, PO Box 88,}\\
{\small         Manchester M60 1QD,}\\
{\small         United Kingdom}\\
{\small         \texttt{richard.booth@umist.ac.uk}}\\
  \and Alexandre V. Borovik\\
{\small Department of Mathematics,}\\
{\small         UMIST, PO Box 88,}\\
{\small         Manchester M60 1QD,}\\
{\small         United Kingdom}\\
{\small         \texttt{borovik@umist.ac.uk}}\\
  \and Neil White\thanks{Supported by EPSRC Grant GR/R53593.}\\
  {\small Department of Mathematics}\\
  {\small University of Florida}\\
  {\small Gainesville, }\\
  {\small Florida 32611, USA}\\
  {\small \texttt{white@math.ufl.edu}}
        }
\newtheorem{lem}{Lemma}
\newtheorem{thm}[lem]{Theorem}
\newtheorem{cor}[lem]{Corollary}

\newcommand{\qed}{\hfill $\diamond $ \\ \vskip 5pt}
\newenvironment{proof}[1][Proof]{\paragraph{#1}}{\qed}
\newcommand{\defterm}[1]{\emph{#1}}

\newcommand{\bigp}[1]{\left( #1 \right)}

\newcommand{\bigcp}[1]{\left\{ #1 \right\}}

\newcommand{\B}{\ensuremath{{\cal B}}}
\newcommand{\M}{\ensuremath{{\cal M}}}

\newcommand\symdiff{\mathbin{\mbox{$\bigtriangleup$}}}
\begin{document}
\maketitle

Represented Coxeter matroids of types $C_n$ and $D_n$, that is,
symplectic and orthogonal matroids arising from totally isotropic
subspaces of symplectic or (even-dimensional) orthogonal spaces, may
also be represented in buildings of type $C_n$ and $D_n$, respectively
(see \cite[Chapter 7]{4}).  Indeed, the particular buildings involved
are those arising from the flags or oriflammes, respectively, of
totally isotropic subspaces.  There are also buildings of type $B_n$
arising from flags of totally isotropic subspaces in odd-dimensional
orthogonal space. Coxeter matroids of type $B_n$ are the same as those
of type $C_n$ (since they depend only upon the directions of roots,
not the length of roots---see the Gelfand-Serganova Theorem below).
However, buildings of type $B_n$ are distinct from those of the other
types. Thus the question arises whether there are such things as
$B_n$-representable matroids, that is, those representable in odd
dimensional orthogonal space, and presumably therefore in such
buildings. We answer this question in the affirmative in \cite{2}. The
matroids so representable turn out to be a special case of symplectic
(flag) matroids, those whose top component, or Lagrangian matroid, is
a union of two Lagrangian orthogonal matroids. These two matroids are
called a Lagrangian pair, and they are the combinatorial manifestation
of the ``fork'' at the top of an oriflamme (or of the fork at the end
of the Coxeter diagram of $D_n$). Thus Lagrangian pairs are a very
natural subject of investigation.

Here we give a number of equivalent characterizations of Lagrangian
pairs, and prove some rather strong properties of them.

\section{Symplectic and Orthogonal Matroids}

Let
$$
[n]=\{1,2,\ldots,n\} \;  \hbox{and}\; [n]^{*}=\{1^{*},2^{*},\ldots,n^{*}\}.
$$
Define the map $*:[n]\rightarrow [n]^{*}$ by $i\mapsto i^{*}$ and
the map $*:[n]^{*}\rightarrow [n]$ by $i^{*}\mapsto i $. In other
words, we are defining $i^{**}=i$. Then $*$ is an involutive
permutation of the set $[n]\cup [n]^*$.

We say that a subset $K\subset [n]\cup[n]^*$ is \emph{admissible} if and only if
 $K\cap K^{*}=\emptyset $.

A linear ordering $\prec$ of $[n]\cup[n]^*$ is called a {\em
$C_n$-admissible ordering} if  $i\prec j$ implies
that\/ $j^* \prec i^*$ for all $i,j\in[n]\cup[n]^*$.
Equivalently, an ordering $\prec$ on $[n]\cup[n]^*$ is
$C_n$-admissible if and only if,
when the $2n$ elements are listed from largest to smallest,
the first $n$ elements listed form an admissible set, and the last $n$
elements listed are the stars of the first $n$ elements listed,
but are listed
in reverse order. A \emph{$D_n$-admissible ordering} of $[n]\cup[n]^*$ is
similar to a $C_n$-admissible ordering, except that the middle two elements
(i.e., the $n$-th and $n+1$-st elements in the above listing) are now
incomparable.

Denote by $J_k$ the collection of all admissible $k$-subsets in $J$,
for some $k\leqslant n$.
If $\prec$ is $C_n$ or $D_n$-admissible ordering on $[n]\cup[n]^*$,
it induces the partial ordering (which we denote by the same symbol
$\prec$) on $J_k$: if $A,B \in J_k$ and
$$
A=\{a_1\prec a_2\prec \cdots \prec a_k\} \quad\hbox{ and }\quad
B=\{b_1\prec b_2\prec \cdots \prec b_k\},
$$
we set $A\prec B$ if
$$
a_1 \prec b_1, a_2 \prec b_2, \ldots, a_k \prec b_k.
$$
This partial ordering is called the Gale ordering on $J_k$ induced by $\prec$.
Given an arbitrary partial order $\leqslant$ on a set $T$,
we can also induce a partial
ordering $\leqslant$ on $T\times T$ by saying
$(s,t)\leqslant (u,v)$ whenever $s\leqslant u$
and $t\leqslant v$. We can likewise induce a
partial order on unordered pairs of
elements of $T$ by setting $\{\,s,t\,\}\leqslant\{\,u,v\,\}$
whenever $(s,t)\leqslant
(u,v)$ or $(s,t)\leqslant (v,u)$.

Now let $\B \subseteq J_k$
be a collection of admissible $k$-element subsets of the set $J$.
We say that $M=(^*,\,\B)$ is a {\em
symplectic matroid\/} if it satisfies the following {\em
Maximality Property}:
\begin{quote}
\emph{for every $C_n$-admissible order $\prec$ on $J$, the collection
$\B$ contains a unique maximal member, i.e.\ a subset $A \in
\B$ such that  $B \prec A$ {\rm (}in the Gale order
induced by $\prec${\rm )}, for all $B\in \B$}.
\end{quote}
The collection $\B$ is called the \emph{ collection of bases} of the symplectic
matroid $M$, its elements are called \emph{bases} of $M$, and the cardinality
$k$ of the bases is the \emph{rank\/} of $M$. An \emph{orthogonal
matroid\/} is defined similarly using $D_n$-admissible orderings.
Ordinary matroids on $[n]$ can be defined in the similar fashion, using
$A_n$-admissible orderings, which are arbitrary linear orderings on $[n]$;
indeed, this is essentially the well-known greedy algorithm of matroid theory.
A \emph{Lagrangian
matroid} (resp. \emph{Lagrangian orthogonal matroid\/}) is a symplectic matroid
(resp. orthogonal matroid\/) of rank $n$, the maximum possible.

A useful characterization of symplectic and orthogonal
matroids is given by the Gelfand-Serganova Theorem (see \cite{4}).
For $B\in\B$, define a point in a real vector
space spanned by $\{\,\epsilon_i: i\in[n]\,\}$ by
$$\delta_B=\sum_{j\in B}\epsilon_j,$$
where $\epsilon _{i^*}$ is
defined to be $-\epsilon_i$. Then $\Delta_\B$ is defined to be
the convex hull of the $\delta_B$ for $B\in\B$.  We define roots
for $D_n$ to be all vectors of the form $\epsilon_j-\epsilon_k$ for
$j,k\in[n]\cup[n]^*$. Roots for $C_n$ are the same together with all
vectors of the form $2\epsilon_j$ for $j\in[n]\cup[n]^*$. Then the
\emph{Gelfand-Serganova Theorem} says that if $\B\subseteq J_k$,
then $\B$ is a symplectic (resp. orthogonal) matroid if and only
if $\Delta_\B$ has all of its edges (i.e., one-dimensional
faces) parallel to roots for $C_n$ (resp. $D_n$).

In the case of a Lagrangian orthogonal matroid, if $\delta_A$ and $\delta_B$
are adjacent vertices in $\Delta$, then the edge between them is parallel to a
root $\epsilon_j-\epsilon_k$. Since $A$ and $B$ are admissible $n$-sets, each
must have either $j$ or $j^*$ as an element, and likewise $k$ or $k^*$.  It
follows that $B=(j,k)(j^*,k^*)A$, regarded as a permutation (in cycle
notation) acting on $A$. Whether $j$ and $k$ are in $[n]$ or $[n]^*$,
it follows that $B$ has the same parity as $A$ in terms of number of starred
elements. Consequently, all bases of a Lagrangian orthogonal matroid have the
same parity. The same is not true for non-Lagrangian orthogonal matroids; for
example $(1,2^*)(1^*,2)\{1\}=\{2^*\}$.

Another characterization of Lagrangian orthogonal matroids comes from
cosets in the group $D_n$. As a permutation group on $[n]\cup[n]^*$, $D_n$ is
generated by the involutions
\begin{eqnarray*}
s_1&=&(1,2)(1^*,2^*) \\
s_2&=&(2,3)(2^*,3^*)\\
&\vdots&\\
s_{n-1}&=&(n-1,n)((n-1)^*,n^*)\\\
s_n&=&(n-1,n^*)((n-1)^*,n)
\end{eqnarray*}
in cycle notation. Consider the two maximal parabolic subgroups $P^{n}$ and
$P^{n-1}$, generated by $\{\,s_1,s_2,\ldots,s_{n-1}\,\}$ and $\{\,s_1,s_2,
\ldots,s_{n-2},s_n\,\}$, respectively. Notice that $P^n$ is the stabilizer in
$D_n$ of the
admissible $n$-set $[n]$, and hence the left cosets of $P^n$ in $D_n$
correspond to the orbit of $[n]$, namely, all admissible $n$-sets of even
parity. Similarly, $P^{n-1}$ is the stabilizer of $\{\,1,2,\ldots,n-1,n^*\,\}$,
and its left cosets correspond to all admissible $n$-sets of odd parity.
Furthermore, in similar fashion, left cosets of $P^n\cap P^{n-1}$ can be shown
to correspond to admissible $(n-1)$-sets.
Letting elements of $D_n$ act on $D_n$-admissible orderings in the obvious way,
we find that $D_n$ corresponds bijectively to the set of all $D_n$-admissible
orderings. Now we can characterize Lagrangian orthogonal matroids $M$ of even
parity as maps from $D_n$ to $D_n/P^n$ (the set of all left cosets of $P_n$ in
$D_n$). An element $\sigma$ of $D_n$ is sent to the coset corresponding to the
maximal basis of $M$ given by the Maximality Property for the admissible order
corresponding to $\sigma$. This \emph{matroid map} $\mu:D_n\rightarrow D_n/P^n$
is actually very natural,
as the Maximality Property can be entirely rephrased in terms of the cosets
using Bruhat order, see \cite{4}. Likewise, Lagrangian orthogonal matroids of odd
parity give matroid maps $\mu:D_n \rightarrow D_n/P^{n-1}$.

One more very useful characterization of Lagrangian orthogonal matroids is the
Strong Exchange Property \cite{1}. A collection $\B\subseteq J_n$ is the
collection of bases of a Lagrangian orthogonal matroid if and only if:
\begin{quote}
For every $A,B\in\B$ and $a\in A\symdiff B$, there exists $b\in B\smallsetminus
A$ with $b\not=a^*$, such that both $A\symdiff\{\,a,b,a^*,b^*\,\}$ and
$B\symdiff\{\,a,b,a^*,b^*\,\}$ are members of $\B$.
\end{quote}

\section{Characterisations of Lagrangian Pairs}

Consider an admissible set of size $n-1$.  Such a set can be completed
to an admissible set of size $n$ in exactly two ways, by appending
either $i$ or $i^*$ for some $i$.  The two resulting sets are called a
\defterm{Lagrangian pair} of sets, and are characterised by the fact
that their symmetric difference is exactly $\bigcp{i,i^*}$.

Consider now two Lagrangian orthogonal matroids $\M_1$, $\M_2$ of rank
$n$ and of opposite parity.  We say that they form a Lagrangian pair
(of Lagrangian orthogonal matroids) if they satisfy:
\begin{quote}
  For every admissible ordering, the maximal bases of $\M_1$ and
  $\M_2$ under the ordering are a Lagrangian pair of sets.
\end{quote}

We say that a pair of Lagrangian subspaces of orthogonal $2n$-space
form a Lagrangian pair of subspaces if their intersection is of
dimension $n-1$. The following result is well-known.

\begin{lem}
  \label{lem:oriflamme}
  A totally isotropic subspace of dimension $n-1$ in orthogonal
  $2n$-space is contained in exactly two Lagrangian subspaces {\rm (}which
  are a Lagrangian pair{\rm )}.
\end{lem}

\begin{thm}
  A Lagrangian pair of subspaces represent a Lagrangian pair of
  orthogonal matroids.
\end{thm}
\begin{proof}
Let $U_1$ and $U_2$ be a Lagrangian pair of subspaces of a  $2n$-dimensional
orthogonal space. Let a $D_n$-admissible ordering be given. The totally
isotropic subspace $U_1\cap U_2$ can be represented by an $(n-1)\times 2n$ matrix
$C$, with columns indexed by elements of $[n]\cup[n]^*$, see \cite[Chapter 3]{4}.
Reorder the columns
so that they are in the given ordering, where the two columns indexed by the
two unrelated elements may be put in either order. Let $C'$ be the reduced row
echelon form of $C$. Then the pivot columns of $C'$ are those indexed by the
maximal basis $B$ of the rank $n-1$ orthogonal matroid represented by $C$, and,
in particular, the pivot columns must be indexed by an admissible set.
Either $U_1$ or $U_2$ may be similarly represented by adding one row to $C$
and again row reducing. However, by elementary linear algebra, the original
$n-1$ pivot columns remain pivot columns in each case, with one additional
pivot column being added in each case. Thus the maximal bases of the Lagrangian
orthogonal matroids represented by $U_1$ and $U_2$ are both admissible $n$-sets
containing $B$, hence are a Lagrangian pair of sets. Thus these two matroids
are a Lagrangian pair.
\end{proof}

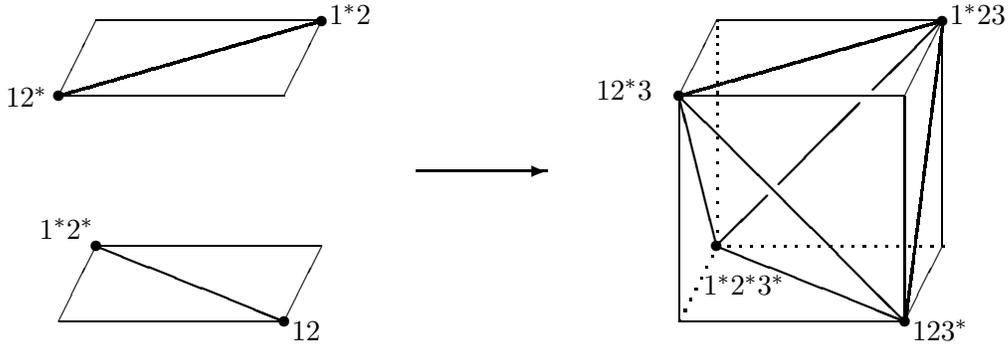
\begin{figure}[h]
\setlength{\unitlength}{0.5mm}
\begin{picture}(100,130)(0,-20)

\thinlines
\put(15,0){\line(1,0){60}}
\put(25,20){\line(1,0){60}}
\put(15,60){\line(1,0){60}}
\put(25,80){\line(1,0){60}}

\put(15,0){\line(1,2){10}}
\put(15,60){\line(1,2){10}}
\put(75,0){\line(1,2){10}}
\put(75,60){\line(1,2){10}}

\thicklines
\qbezier(15,60)(50,70)(85,80)
\put(25,20){\line(5,-2){50}}

\put(110,40){\vector(1,0){35}}

\put(15,60){\circle*{3}}
\put(85,80){\circle*{3}}
\put(75,0){\circle*{3}}
\put(25,20){\circle*{3}}

\put(1,57){\small $12^*$}
\put(10,22){\small $1^*2^*$}
\put(87,79){\small $1^*2$}
\put(77,-5){\small $12$}

\thinlines
\put(180,0){\line(1,0){60}}
\put(180,60){\line(1,0){60}}
\put(190,80){\line(1,0){60}}

\put(180,0){\line(0,1){60}}
\put(240,0){\line(0,1){60}}
\put(250,20){\line(0,1){60}}

\put(180,60){\line(1,2){10}}
\put(240,0){\line(1,2){10}}
\put(240,60){\line(1,2){10}}

\thicklines
\qbezier(180,60)(215,70)(250,80)
\qbezier(240,0)(245,40)(250,80)
\put(180,60){\line(1,-1){60}}
\put(190,20){\line(1,1){13.5}}
\put(206.5,36.5){\line(1,1){21.5}}

\put(250,80){\line(-1,-1){18.5}}
\put(190,20){\line(-1,4){10}}
\put(190,20){\line(5,-2){50}}

\qbezier[20](190,20)(220,20)(250,20)
\qbezier[20](190,20)(190,60)(190,80)
\qbezier[7](180,0)(185,10)(190,20)

\put(180,60){\circle*{3}}
\put(250,80){\circle*{3}}
\put(240,0){\circle*{3}}
\put(190,20){\circle*{3}}

\put(158,59){\small $12^*3$}
\put(186.5,7){\small $1^*2^*3^*$}
\put(252,79){\small $1^*23$}
\put(242,-5){\small $123^*$}

\end{picture}
\caption[Exploded sum of a Lagrangian pair]{{\small\sl The exploded sum of a Lagrangian
pair of Lagrangian orthogonal matroids.}}
\label{fig:expl-sum}
\end{figure}

Let $M_1,M_2$ be Lagrangian orthogonal matroids on $[n]\cup[n]^*$ of
opposite parity, $\B_1,\B_2$ their collections of bases,
and $$\mu_1:W\rightarrow D_n/P^n,\quad \mu_2:W\rightarrow D_n/P^{n-1}$$ the corresponding
matroid maps.
Let
\begin{eqnarray*}
\B_1+(n+1)&=&\{\,B\cup\{\,n+1\,\}\mid B\in\B_1\,\},\\
\B_2+(n+1)^*&=&\{\,B\cup\{\,(n+1)^*\,\}\mid B\in\B_2\,\},\\
\B_3&=&\bigp{\B_1+(n+1)}\cup\bigp{\B_2(n+1)^*},
\end{eqnarray*}
and
\begin{eqnarray*}
\B_4&=&\{\,A\mid A\subseteq [n]\cup[n]^*,\;|A|=n-1,\\
&& \qquad \qquad \hbox{\rm and there exist }
B_1\in\B_1, B_2\in\B_2\\
&&  \qquad \qquad\qquad \qquad\hbox{\rm such that}\quad A=B_1\cap B_2\,\}.
\end{eqnarray*}
We call $\B_3$ the (collection of bases of the)
\emph{exploded sum} of $M_1$ and $M_2$, see Figure~\ref{fig:expl-sum}.

\begin{thm}
\label{thm:tfae}
The following are equivalent:
\begin{enumerate}
\item[{\rm (1)}] $M_1$ and $M_2$ are a Lagrangian pair,
\item[{\rm (2)}] for all $w\in W$, $\mu_1(w)\cap\mu_2(w)\not=\emptyset,$
\item[{\rm (3)}] $\B_3$ is the collection of bases of a
  Lagrangian orthogonal matroid,
\item[{\rm (4)}] $\B_4$ is the collection of bases of a
  Lagrangian orthogonal matroid, and for each $B\in \B_1\cup \B_2$ there
  exists $X\in \B_4$ with $X\subset B$.
\item[{\rm (5)}] for all $w\in W$, there exists a unique unordered
  pair $\{\,B_1,B_2\,\}$, with $B_1\in\B_1, B_2\in\B_2$,
  such that for all $A_1\in\B_1, A_2\in\B_2$, we have
  $\{\,B_1,B_2\,\}\leqslant^w\{\,A_1,A_2\,\},$
\item[{\rm (6)}] for all $w\in W$, there exists a unique ordered pair
  $(B_1,B_2)$, with $B_1\in\B_1, B_2\in\B_2$, such that
  for all $A_1\in\B_1, A_2\in\B_2$, we have
  $(B_1,B_2)\leqslant^w(A_1,A_2).$
\end{enumerate}
\end{thm}
\begin{proof}
The equivalence (4) $\Leftrightarrow$ (5) is proved in [6]. Furthermore,
(1) $\Leftrightarrow$ (6), (6) $\Rightarrow$ (5), and (5) $\Rightarrow$ (2)
are immediate. We will now prove (1) $\Rightarrow$ (3), (3) $\Rightarrow$ (1),
and (2) $\Rightarrow$ (1).

Assume that $M_1$ and $M_2$ are a Lagrangian pair, and choose a
$D_{n+1}$-admissible order on $[n+1]\cup[n+1]^*$.  Restrict this ordering to
$[n]\cup[n]^*$. Notice that this may now be either a $C_n$ or $D_n$-admissible
order. Let $B_1\in\B_1$, $B_2\in\B_2$ be the maximal bases in this
restricted order. If the restricted order is a $C_n$-admissible order, it can
be changed to a $D_n$-admissible order by deleting the relation between the
pair of elements in the middle, and clearly the maximal bases remain unchanged.
Hence $B_1\symdiff B_2=\{\,i,i^*\,\}$.  Clearly $B_1'=B_1\cup\{\, n+1\,\}$ and
$B_2'=B_2\cup\{\,(n+1)^*\,\}$ are the only two candidates for maximal members
of $\B_3$. But $B_1'\symdiff B_2'=\{\,i,i^*,n+1,(n+1)^*\,\}$, and in all
possible $D_{n+1}$-admissible orders, $B_1'$ and $B_2'$ are related. Thus
$\B_3$ has a unique maximal member, proving (3).

Now assume (3), and let a $D_n$-admissible order be given.  Let $B_1$ (resp.
$B_2$) be the maximal basis in $\B_1$ (resp. $\B_2$). Extend the
given order arbitrarily to a $D_{n+1}$-admissible order, and let $f$ be a
linear functional compatible with the extended order. By this we mean that
$f$ is a linear functional on the real vector space spanned by the
basis $\{\,\epsilon_1,\epsilon_2,\ldots,\epsilon_{n+1}\,\}$, with $\epsilon_i=
-\epsilon_{i^*}$ by definition, such that $i\prec j$ in the extended order
implies $f(i)<f(j)$, for all $i,j\in[n+1]\cup[n+1]^*$. Now $B_1'=B_1\cup
\{\,n+1\,\}$ (resp. $B_2'=B_2\cup\{\,(n+1)^*\,\}$) is clearly the maximal
basis in $$\B_1'=\{\,B\cup\{\,n+1\,\}\mid B\in\B_1\,\}$$ (resp.
$\B_2'=\{B\cup\{\,(n+1)^*\,\}\mid B\in\B_2\,\}$). Thus $\delta_{B_1'}$
(resp. $\delta_{B_2'}$) is the unique maximal vertex of $\Delta_{\B_1'}$
(resp. $\Delta_{\B_2'}$) under $f$. Since $\Delta_{\B_1'}$ and
$\Delta_{\B_2'}$ lie in parallel hyperplanes, it is easy to see that
$\delta_{B_1'}\delta_{B_2'}$ must be an edge of $\Delta_{\B_3}$. By the
Gelfand-Serganova Theorem, $B_1'$ and $B_2'$ must be related by an exchange
of the form $(n+1,(n+1)^*)(i,i^*)$ (in cycle notation), for some $i\leqslant n$.
It follows that $B_1$ and $B_2$ must be related by the exchange $(i,i^*)$,
proving that $\B_1$ and $\B_2$ are a Lagrangian pair.

Finally, assume (2).  Thus $\mu_1(w)=aP^n$ and $\mu_2(w)=bP^{n-1}$ with
$aP^n\cap bP^{n-1}\not=\emptyset$. But this means $aP^n\cap bP^{n-1}=c(P^n\cap
P^{n-1})=cP^{n,n-1}$. It follows that there is an admissible $(n-1)$-set contained
in both $B_1=\max^w\B_1$ and $B_2=\max^w\B_2$, and thus that $B_1$
and $B_2$ differ by an exchange of the form $(i,i^*)$ for some $i$, showing
that $M_1$ and $M_2$ are a Lagrangian pair.
\end{proof}

Condition (2) in the previous theorem amounts to saying that the two Lagrangian
orthogonal matroids are \emph{concordant}, see [3].

\section{Further results on Lagrangian pairs}

We now need to recall the concept of quotient (or, essentially, strong map)
of ordinary matroids. If $M_1$ and $M_2$ are matroids on the same set $[n]$,
then we say that $M_2$ is a quotient of $M_1$ if every circuit of $M_1$
is a union of circuits of $M_2$. As is shown in [4, Chapter 1], $M_2$ is a
quotient of $M_1$ if and only if, for every linear ordering of $[n]$,
the maximal basis of $M_2$ is a subset of the maximal basis of $M_1$.
To relate ordinary matroids to Lagrangian orthogonal matroids, we need a
mapping $\Phi$ defined as follows: For $B\subseteq[n]$, let $\Phi(B)=
B\cup([n]\smallsetminus B)^*$. Then if $\B$ is the collection of bases of a
matroid, $\Phi(\B)=\{\,\Phi(B)\mid B\in\B\,\}$ is a Lagrangian
orthogonal matroid, as proved in [4, Chapter 3].

\begin{thm}
Let $M_1$ and $M_2$ be ordinary matroids on $[n]$, of ranks $k$ and $k-1$
respectively. Then $M_2$ is a quotient of $M_1$ if and only if $\Phi(M_1)$ and
$\Phi(M_2)$ are a Lagrangian pair.
\end{thm}
\begin{proof}
Let $\prec$ be a $D_n$-admissible ordering of $[n]\cup[n]^*$, and let $\leqslant$
denote the restriction of this ordering to $[n]$, which must be a linear order.
First we claim that if $A$ and $B$ are bases of the same matroid $M$, then
$A\leqslant B$ implies $\Phi(A)\prec\Phi(B)$.  Indeed, if $A\leqslant B$, then
$[n]\smallsetminus A\geqslant [n]\smallsetminus B$, and hence $([n]\smallsetminus A)
^*\leqslant([n]\smallsetminus B)^*$, so $\Phi(A)\prec \Phi(B)$. It follows that $B$
is the maximal basis of $M$ if and only if $\Phi(B)$ is the maximal basis of
$\Phi(M)$.

Now let $B_1$ and $B_2$ be the maximal bases of $M_1$ and $M_2$, resp.  If
$M_2$ is a quotient of $M_1$, then $B_2\subseteq B_1$, say $B_2=B_1
\smallsetminus\{\,i\,\}$. Then $\Phi(B_1)$ and $\Phi(B_2)$ differ by the
exchange $(i,i^*)$.
But we have just seen that these are the maximal bases of $\Phi(M_1)$ and
$\Phi(M_2)$, resp., showing that these two Lagrangian orthogonal matroids are
a Lagrangian pair.

Conversely, suppose that $M_2$ is not a quotient of $M_1$. Then there exists a
linear ordering on $[n]$ such the maximal bases $B_1$ and $B_2$ of $M_1$ and
$M_2$ (resp.) are not related by containment. Extend this ordering to a
$D_n$-admissible ordering on $[n]\cup[n]^*$.  Then the maximal bases
$\Phi(B_1)$ and $\Phi(B_2)$ of $\Phi(M_1)$ and $\Phi(M_2)$ (resp.) are not
related by an exchange of the form $(i,i^*)$. Thus $\Phi(M_1)$ and $\Phi(M_2)$
are not a Lagrangian pair.
\end{proof}

As a corollary, we get a new characterization of elementary quotients.

\begin{cor}
Let\/ $\B_1$ and\/ $\B_2$ be the collection of bases of two matroids
$M_1$ and $M_2$ of ranks $k$ and\/ $k-1$, respectively.
Then $M_2$ is an elementary
quotient of\/ $M_1$ if and only if for every $B_1\in\B_1, B_2\in
\B_2,$ and\/ $i\in B_1\symdiff B_2$, either
\begin{itemize}
\item[{\rm (1)}] $B_1\symdiff\{\,i\,\}\in \B_2, B_2\symdiff\{\,i\,\}\in \B_1$, or
\item[{\rm (2)}] there exists $j\in B_1\symdiff B_2$ such that
$$B_1\symdiff\{\,i,j\,\}\in
\B_1 \quad\hbox{ and }\quad B_2\symdiff\{\,i,j\,\}\in\B_2.$$
\end{itemize}
Furthermore, for every $B_1\in\B_1, B_2\in\B_2$ there exists $i\in
B_1\symdiff B_2$ such that case {\rm (1)} holds.
\end{cor}

\begin{proof}
We have that $M_2$ is an elementary quotient of $M_1$ if and only if
$\Phi(M_1)$ and $\Phi(M_2)$ are a Lagrangian pair if and only if their
exploded sum is a Lagrangian orthogonal matroid, which is true if and only if
their exploded sum satisfies the Strong Exchange Property.  The
corollary follows immediately from translating what the Strong Exchange
Property says in terms of $\B_1$ and $\B_2$.
\end{proof}

\begin{thm}
\label{thm:proj}
Let\/ $\B$ be a Lagrangian orthogonal matroid on $[n]\cup[n]^*$, and\/ $i\in[n]$.
Define
$$\B_1=\{\,B\smallsetminus\{\,i\,\}\mid i\in B, B\in\B\,\},$$
$$\B_2=\{\,B\smallsetminus\{\,i^*\,\}\mid i^*\in B, B\in\B\,\}.$$
Then $\B_1$, $\B_2$ are a Lagrangian pair of Lagrangian orthogonal matroids on
$([n]\cup[n]^*)\smallsetminus\{\,i,i^*\,\}$.
\end{thm}

\begin{proof}
We see that $\B$ is just the exploded sum of $\B_1$ and $\B_2$, after
appropriate relabelling, so the desired result follows from
Theorem~\ref{thm:tfae}.
\end{proof}

\begin{thm}
  Let $\B_1$, $\B_2$ be the collections of bases of a Lagrangian pair of
  Lagrangian orthogonal matroids.  Then $\B=\B_1\cup\B_2$ is the collection of
  bases of
  a Lagrangian (symplectic) matroid.
\end{thm}
\begin{proof}
If $\B_1$, $\B_2$ are a Lagrangian pair, then their exploded sum $\B_3$ is a
Lagrangian orthogonal matroid, by Theorem~\ref{thm:tfae}.
The Strong Exchange Axiom on $\B_3$ now implies the Symmetric Exchange
Axiom on $\B_1\cup \B_2$.
\end{proof}

Since the union of a Lagrangian pair is a symplectic matroid,
it is natural to wonder whether, given Lagrangian orthogonal matroids
of the same rank and opposite parity whose union is a symplectic
matroid, they are necessarily a Lagrangian pair.  Figure~\ref{fig:notlp}
shows that the answer is no. Indeed, since $123$ is the only basis of the
matroid of even parity, any admissible order which makes $1^*2^*3^*$ the
maximal basis of the other matroid violates the definition of Lagrangian pair.

\unitlength 1mm
\begin{figure}
\begin{picture}(80,80)(-20,0)
{\tiny
\put(51,9){\makebox(0,0)[tl]{$123^*$}}
\put(9,51){\makebox(0,0)[br]{$12^*3$}}
\put(29,25){\makebox(0,0)[r]{$1^*2^*3^*$}}
\put(49,51){\makebox(0,0)[br]{$123$}}
\put(71,66){$1^*23$}
}
\thicklines
\put(10,50){\line(1,0){39}}
\put(50,10){\line(0,1){39}}
\qbezier[20](10,10)(30,10)(50,10)
\qbezier[20](10,10)(10,30)(10,50)
\qbezier[20](30,65)(50,65)(70,65)
\qbezier[20](70,25)(70,45)(70,65)
\qbezier[20](30,25)(50,25)(70,25)
\qbezier[20](30,25)(30,45)(30,65)
\put(50.75,50.78){\line(4,3){19}}
\qbezier[15](10,50)(20,57.5)(30,65)
\qbezier[15](10,10)(20,17.5)(30,25)
\qbezier[15](50,10)(60,17.5)(70,25)
\thicklines
\put(30,25){\line(1,1){1.5}}
\put(33.5,28.5){\line(1,1){36.5}}
\put(30,25){\line(4,-3){20}}
\put(30,25){\line(-4,5){20}}
\put(10,50){\line(1,-1){40}}
\put(10,50){\line(4,1){60}}
\qbezier(50,10)(60,37.5)(70,65)
\put(50,10){\circle*{2.5}}
\put(10,50){\circle*{2.5}}
\put(30,25){\circle*{2.5}}
\put(70,65){\circle*{2.5}}
\put(50,50){\circle{2.5}}
\thinlines
\end{picture}
\caption{{\small\sl A Lagrangian symplectic matroid which is not the union of a
Lagrangian pair.}}
\label{fig:notlp}
\end{figure}
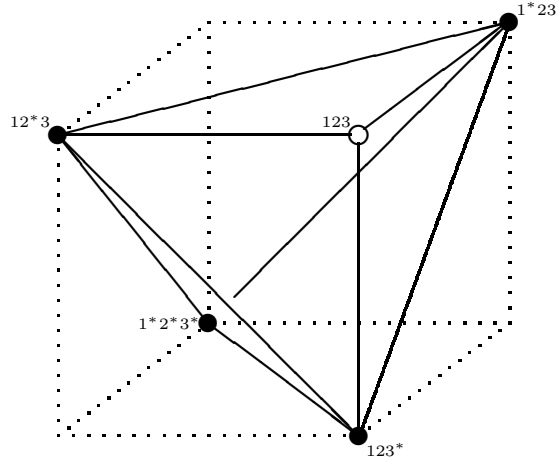

\begin{thm}
Let $\B_1$ be the collection of bases of a Lagrangian orthogonal matroid,
and $\B_2=(i,i^*)\B_1$, where $(i,i^*)$ is a transposition,
to be thought of as a permutation acting elementwise on the members of $\B_1$.
Then $\B_1$ and $\B_2$ are a Lagrangian pair.
\end{thm}

\begin{proof}
By Theorem~\ref{thm:tfae}, it suffices to prove that the exploded sum,
$$\B_3=\B_1+(n+1)\cup\B_2+(n+1)^*$$
is a Lagrangian orthogonal matroid. So we wish to show that $\B_3$
satisfies the Strong Exchange Property. Let $A,B\in\B_3$ and
$a\in A\smallsetminus B$. We need to show that there exists $b\in (A\symdiff B)
\smallsetminus\{\,a,a^*\,\}$ so that $A\symdiff\{\,a,a^*,b,b^*\,\}$
and $B\symdiff
\{\,a,a^*,b,b^*\,\}$ are both in $\B_3$. If both $A$ and $B$ are in
$\B_1+(n+1)$ or else both in $\B_2+(n+1)^*$, then we are done, for
$a\not=(n+1),(n+1)^*$, and Strong Exchange holds in $\B_1$ and in
$\B_2$. Thus we may assume that $A\in\B_1+(n+1)$ and $B\in\B_2
+(n+1)^*$.

Case 1. Suppose $a=i$ or $i^*$.  Then we may choose $b=n+1$, and we are done.

Case 2. Suppose $a=n+1$.  If $i\in A\symdiff B$, then we choose $b=i$ and we are
done. Thus we suppose, without loss of generality, that $i\in A\cap B$.  Let
$B'=B\symdiff\{\,i,i^*,n+1,(n+1)^*\,\}$.  Then $B'\in\B_1+(n+1)$, and $i\in
A\symdiff B'$. By Strong Exchange on $\B_1$, we have $j\in [n]\cup[n]^*$
so that $$A':=A\symdiff \{\,i,i^*,j,j^*\,\}$$ and $$B'\symdiff\{\,i,i^*,j,j^*\,\}=
B\symdiff\{\,j,j^*,n+1,(n+1)^*\,\}$$ are both in $\B_1+(n+1)$.  Thus
\begin{eqnarray*}
A\symdiff \{\,j,j^*,n+1,(n+1)^*\,\} &=&
A'\symdiff\{\,i,i^*,n+1,(n+1)^*\,\}\\
&\in& \B_2+ (n+1)^*.
\end{eqnarray*}
Thus $b=j$ gives the desired Strong Exchange.

Case 3. We are left with $a\not=i,i^*,n+1,(n+1)^*$. Let
$$B'=b\symdiff\{\,i,i^*,
n+1,(n+1)^*\,\}\in \B_1+(n+1).$$ Since $a\in A\symdiff B'$, there exists
$j\in A\symdiff B'$, $j\not=a,a^*,n+1,(n+1)^*,$ so that
$$A':=A\symdiff\{\,a,a^*,j,j^*
\,\}\in \B_1+(n+1)$$ and $B'\symdiff\{\,a,a^*,j,j^*\,\}\in\B_1+(n+1)$.
Then
\begin{eqnarray*}
B\symdiff\{\,a,a^*,j,j^*\,\}&=& B'\symdiff\{\,i,i^*,n+1,(n+1)^*\,\}\symdiff
\{\,a,a^*,j,j^*\,\}\\
&\in& \B_2+(n+1)^*,
\end{eqnarray*} regardless of whether $j=i$ or $i^*$, or not.
Thus $b=j$ again gives the desired Strong Exchange.
\end{proof}



\end{document}